\documentclass[a4paper] {article} [12pt]
\usepackage{amsmath,amsfonts,amssymb,amsthm,amscd,nccmath}
\usepackage{graphicx}
\usepackage[numbers]{natbib}
\usepackage{color}

\newtheorem{Theorem}{Theorem}[section]
\newtheorem{thm}[Theorem]{Theorem}
\newtheorem{Proposition}[Theorem]{Proposition}
 
\newtheorem{defn}[Theorem]{Definition}

\newtheorem{conj}[Theorem]{Conjecture}

\newtheorem{lem}[Theorem]{Lemma}

\newtheorem{Fact}[Theorem]{Fact}

\makeatletter
\newcommand{\thorn}{{\fontencoding{T1}\selectfont\th}}
\newcommand{\@indepsymbol}[2]{#1\setbox0=\hbox{$#1x$}\kern\wd0\hbox to 0pt{\hss$#1\mid$\hss}\lower.9\ht0\hbox to 0pt{\hss$#1\smile$\hss}\kern\wd0}
\newcommand{\@nindepsymbol}[2]{#1\setbox0=\hbox{$#1x$}\kern\wd0\hbox to 0pt{\mathchardef
	\nn=12854\hss$#1\nn$\kern1.4\wd0\hss}\hbox to
	0pt{\hss$#1\mid$\hss}\lower.9\ht0 \hbox to
	0pt{\hss$#1\smile$\hss}\kern\wd0}
\newcommand{\ind}[1][]{\mathop{\mathpalette\@indepsymbol{}^{\!\!\!\!\rlap{$\scriptstyle\textnormal{#1}$}\,\,\,\,}}}
\newcommand{\nind}[1][]{\mathop{\mathpalette\@nindepsymbol{}^{\!\!\!\rlap{$\scriptstyle\textnormal{#1}$}\,\,\,}}}
\newcommand{\@Ind}[1][]{\mathpalette\@indepsymbol{}^{\!\!\!\!\mbox{$\scriptstyle\textnormal{#1}$}}}
\newcommand{\Ind}[1][]{\@Ind[\ \,]}
\newcommand{\newind}[4]{
	\newcommand{#1}{{\!\@Ind[#4]}}
	\newcommand{#2}{\ind[#4]}
	\newcommand{#3}{\nind[#4]}
}
\newind{\thInd}{\thind}{\nthind}{\thorn}
\renewcommand{\thInd}{\text{$\@Ind[\thorn]$\;}}
\makeatother

\title{The definable $(p,q)$-theorem for distal theories}

\date{\today}
\author{Gareth Boxall, Charlotte Kestner}

\begin{document}

\maketitle

\begin{abstract}
Answering a special case of a question of Chernikov and Simon, we show that any non-dividing formula over a model $M$ in a distal NIP theory is a member of a consistent definable family, definable over $M$.  
\end{abstract}

\section{Introduction}

The behaviour of forking in NIP structures is a mysterious topic currently under investigation. It is well-known that forking provides a good notion of independence for stable and even simple structures and that it fails to do so in the unstable NIP setting in some quite significant ways. However, certain aspects of the good behaviour of forking in stable structures do extend to NIP and there has been considerable work in identifying which these are. A notable example is the equivalence of forking and dividing over models, which Chernikov and Kaplan showed to be true even in the more general NTP$_2$ setting in \cite{CK}. Conjecture \ref{mainconjecture}, if true, would be another interesting example. It is stated in \cite{CS} (as a question, but has subsequently achieved the status of a conjecture). As is discussed in Section 2 of \cite{Sinv}, it can be seen as a definable analogue of the $(p,q)$-theorem from combinatorics.

\begin{conj}\label{mainconjecture}
Let $M\prec N$ be NIP $L$-structures. Let $\varphi(x,y)$ be an $L_M$-formula and let $b\in N^{|y|}$. Assume $\varphi(x,b)$ does not divide over $M$. Then there is an $L_M$-formula $\psi(y)\in tp(b/M)$ such that $\{\varphi(x,b^\prime):b^\prime\in \psi(M)\}$ is consistent. 
\end{conj}

This has been proved in certain special cases by Simon and Starchenko. In \cite{Sdp} Simon proves it under the additional assumption that the theory of $M$ is dp-minimal and has small or medium directionality. In \cite{SS} Simon and Starchenko obtain a stronger conclusion for a certain class of dp-minimal structures including all those with definable Skolem functions. In \cite{Sinv} Simon proves Conjecture \ref{mainconjecture} under the additional assumption that there is some countable $M^\prime\equiv M$ with the property that, for any complete type $p$ over $M^\prime$ and any elementary extension $N^\prime$, only countably many complete types over $N^\prime$ are coheirs of $p$. Conjecture \ref{mainconjecture} is known under the additional assumption of stability. There are also two interesting approximations to Conjecture \ref{mainconjecture} which we mention in \S 2. 

Our main result is to establish Conjecture \ref{mainconjecture} under the additional assumption that $M$ is distal. 

\begin{thm}\label{mainresult}
Let $M\prec N$ be distal NIP $L$-structures. Let $\varphi(x,y)$ be an $L_M$-formula and let $b\in N^{|y|}$. Assume $\varphi(x,b)$ does not divide over $M$. Then there is an $L_M$-formula $\psi(y)\in tp(b/M)$ such that $\{\varphi(x,b^\prime):b^\prime\in \psi(M)\}$ is consistent. 
\end{thm}

We recall a definition of distality in \S 2. Simon introduced this concept in \cite{Sdist} to single out the class of (in some sense) completely non-stable NIP structures. It is encouraging to have the conclusion of Conjecture \ref{mainconjecture} now for both stable and distal NIP structures. With both extremes covered, perhaps the remaining inbetween cases will follow soon. Simon has, in recent work, explored ways in which a type in an NIP structure can be ``decomposed" into a stable part and a distal part. He suggested the idea of using this work to remove the distal assumption in Theorem \ref{mainresult}. This is an appealing idea, but we have not managed to implement it. 

So-called strict Morley sequences play an important role in our proof of Theorem \ref{mainresult} and we have a lemma concerning them which may be of independent interest. We discuss strict Morley sequences in \S 2 and then present our lemma in \S 3. In \S 4 we combine the various ingredients to complete the proof of Theorem \ref{mainresult}. 

We would like to thank Pierre Simon for valuable discussions concerning this work, including pointing out errors in earlier attempts. We would also like to thank the anonymous referee for, among other things, improving the concluding argument and suggesting the present title of the paper. 

\section{Preliminaries}

In this section we recall some background material. We begin with two approximations to Conjecture \ref{mainconjecture}. The first replaces $\psi(y)$ with $tp(b/M)$ and was established in \cite{CK}. The second replaces ``consistent" with something weaker and is a corollary of Proposition 25 in \cite{CS}, in combination with the sentence preceding the statement of that result in \cite{CS}. Both are discussed in \cite{Sinv}. 

\begin{Proposition}\label{firstapprox}
Let $M\prec N$ be NIP $L$-structures. Let $\varphi(x,y)$ be an $L_M$-formula and let $b\in N^{|y|}$. Assume $\varphi(x,b)$ does not divide over $M$. Let $q(y)=tp(b/M)$. Then $\{\varphi(x,b^\prime):N\models q(b^\prime)\}$ is consistent. 
\end{Proposition} 

\begin{Proposition}\label{secondapprox}
Let $M\prec N$ be NIP $L$-structures with $N$ sufficiently saturated. Let $\varphi(x,y)$ be an $L_M$-formula and let $b\in N^{|y|}$. Assume $\varphi(x,b)$ does not divide over $M$. Then there exist an $L_M$-formula $\psi(y)\in tp(b/M)$ and a finite $A_\psi\subseteq N^{|x|}$ such that, for each $b^\prime\in\psi(M)$, there is some $a\in A_\psi$ such that $N\models \varphi(a,b^\prime)$. 
\end{Proposition}

The concept of distality was introduced by Simon in \cite{Sdist} to single out those NIP structures which are, in some sense, completely unstable. There are various equivalent definitions. The following is (almost the same as and clearly equivalent to) the one used in \cite{CS} together with the assumption of NIP (as we shall not want it without that assumption).

\begin{defn}
A structure $M$ is said to be distal if it is NIP and, for all $M\prec N$, indiscernible $(b_i)_{i\in\mathbb{Z}}$ in $N^{|b_0|}$ and $A\subseteq N$, if $(...,b_{-2},b_{-1},b_{1},b_{2},...)$ is indiscernible over $A$ then also  $(...,b_{-2},b_{-1},b_0,b_{1},b_{2},...)$ is indiscernible over $A$. 
\end{defn}

Like NIP, this is really a property of the complete theory of $M$. Theorem 21 in \cite{CS} gives a characterisation of distality for NIP structures. The following result is a special case of one direction of this characterisation (modulo the obvious and well-known fact that an expansion of a distal structure by constants is again distal).

\begin{Proposition}\label{distality}
Let $N$ be a distal $L$-structure. Then, for every $L_N$-formula $\varphi(x,y)$, there exist an $L_N$-formula $\theta(x, z)$ and a natural number $k$ such that $|z|=k|y|$ and, for all $B, C\subseteq N^{|y|}$ and $a\in N^{|x|}$, if 

\begin{enumerate}

\item $|C|\geq 2$,

\item $B\subseteq C$,

\item $B$ is finite and

\item $N\models \varphi(a,b)$ for all $b\in B$
\end{enumerate}

then there is some $c\in C^k$ such that $N\models \theta(a,c)$ and, for all $b\in B$, $N\models \forall x(\theta(x,c)\rightarrow \varphi(x,b))$.
\end{Proposition}

The notion of a strict Morley sequence comes from work of Shelah and plays an important role in \cite{Sinv}. An indiscernible sequence of realisations of a complete type over $M$ is a strict Morley sequence if it is strictly non-forking over $M$. By the fact, from \cite{Shelah}, that non-forking implies invariance (for types over models), one obtains the following statement which we give in place of the original definition.  

\begin{Fact}\label{strictMorley}
Let $M\prec N$ be NIP structures with $N$ sufficiently saturated. Let $q(y)$ be a complete type over $M$. Let $(b_n)_{n\in\mathbb{N}}$ be a strict Morley sequence in $N^{|y|}$ for $q(y)$. Then $(b_n)_{n\in\mathbb{N}}$ is indiscernible over $M$, $b_0\models q(y)$ and, for each $n\in\mathbb{N}$, the following two conditions are satisfied:

\begin{enumerate}

\item $tp(b_{n+1}/Mb_0...b_n)$ extends to a complete type over $N$ which is $M$-invariant and 

\item $tp(b_0...b_n/Mb_{n+1})$ extends to a complete type over $N$ which is $M$-invariant.

\end{enumerate}
\end{Fact}

If $(b_n)_{n\in\mathbb{N}}$ satisfies the conclusion of Fact \ref{strictMorley} then it is clear that any other realisation of $tp(b_0b_1b_2b_3.../M)$ in $N$ will too. The following was established, in greater generality, by Chernikov and Kaplan in \cite{CK}. See also \cite{Sbook}. 

\begin{Fact}\label{existstrictmorley}
Let $M\prec N$ be NIP structures with $N$ sufficiently saturated. Let $q(y)$ be a complete type over $M$. Then there exists a strict Morley sequence $(b_n)_{n\in\mathbb{N}}$ in $N^{|y|}$ for $q(y)$. 
\end{Fact}

\section{Fitting indiscernible sequences around sets}

In this section we note the following lemma which might be of independent interest. It is an almost immediate consequence of Fact \ref{strictMorley}. 

\begin{lem}\label{fitting}
Let $M\prec N$ be NIP structures with $N$ sufficiently saturated. Let $q(y)$ be a complete type over $M$. Let $(b_n)_{n\in\mathbb{N}}$ be a strict Morley sequence for $q(y)$ in $N^{|y|}$. Let $B\subseteq N^{|y|}$ be a finite set of realisations of $q$. Then there is a sequence $(d_n)_{n\in\mathbb{Z}}$ in $N^{|y|}$ such that, for each $d\in B$, $(...,d_{-2},d_{-1},d,d_1,d_2,...)$ is indiscernible over $M$ and has the same EM-types as $(b_n)_{n\in\mathbb{N}}$ over $M$.
\end{lem}

\proof For each positive $k\in\mathbb{N}$ we can use part 2 of Fact \ref{strictMorley} to obtain $d_0,...,d_{k-1}$ and then part 1 of Fact \ref{strictMorley} to obtain $d_{k+1},d_{k+2},...$ such that, for all $d\in B$, $(d_0,...,d_{k-1},d,d_{k+1},...)$ has the same type as $(b_n)_{n\in\mathbb{N}}$ over $M$. The result then follows by compactness.  \endproof

\section{Proof of the main result}

We now combine the various ingredients to give a proof of Theorem \ref{mainresult}. Let $M\prec N$ be distal $L$-structures. Let $\varphi(x,y)$ be an $L_M$-formula and let $b\in N^{|y|}$. Assume $\varphi(x,b)$ does not divide over $M$. Let $q(y)=tp(b/M)$. The conclusion of Theorem \ref{mainresult} follows trivially when $b\in M^{|y|}$. We assume $b\notin M^{|y|}$. So $q$ has infinitely many realisations.  We may assume $N$ is sufficiently saturated.

\begin{Proposition}\label{refinement}
There exist an $L_M$-formula $\theta(x,z)$, a natural number $k$ such that $|z|=k|y|$ and a complete type $r(z)$ over $M$ such that, for each finite set $B\subseteq N^{|y|}$ of realisations of $q(y)$, there is some $c\in N^{|z|}$ such that the following conditions are satisfied.

\begin{enumerate}

\item $N\models r(c)$,

\item $\theta(x,c)$ does not divide over $M$ and

\item for all $d\in B$, $N\models \forall x(\theta(x,c)\rightarrow \varphi(x,d))$.

\end{enumerate}

\end{Proposition}

\proof Let $\theta(x,z)$ and $k$ be as in Proposition \ref{distality}. Let $(b_n)_{n\in\mathbb{N}}$ be a strict Morley sequence for $q(y)$ in $N^{|y|}$. Let $a\in N^{|x|}$ be such that $N\models \varphi(a,b_n)$ for all $n\in\{1,...,k+1\}$. Such $a$ exists because $\varphi(x,b)$ does not divide over $M$. Then there is a tuple $i_1...i_k$ from $\{1,...,k+1\}$ such that $N\models \theta(a,b_{i_1}...b_{i_k})$ and, for all $a^\prime\in N^{|x|}$, if $N\models \theta(a^\prime,b_{i_1}...b_{i_k})$ then, for all $n\in \{1,...,k+1\}$, $N\models \varphi(a^\prime,b_n)$. 

Since forking equals dividing over models, by \cite{CK}, we may assume there is some small $N^\prime$ such that $M\prec N^\prime\prec N$ and, for all $n\in\mathbb{N}$, $b_n\in N^{\prime |y|}$ and $a$ realises a complete type over $N^\prime$ which does not fork over $M$. It then follows that $\theta(x,b_{i_1}...b_{i_k})$ does not divide over $M$. 

Let $i\in\{1,...,k+1\}\setminus\{i_1,...,i_k\}$. Let $B\subseteq N^{|y|}$ be a finite set of realisations of $tp(b/M)$. Let $(d_n)_{n\in\mathbb{Z}}$ in $N^{|y|}$ be indiscernible over $M$ with the same $EM$-types over $M$ as $(b_n)_{n\in\mathbb{N}}$ and with the property that, for all $d\in B$, $(..., d_{-2},d_{-1},d,d_1, d_2,...)$ is indiscernible over $M$. The existence of such a sequence comes from Lemma \ref{fitting}. Taking $r=tp(b_{i_1}...b_{i_k}/M)$ and then $c=d_{(i_1-i)}...d_{(i_k-i)}$ yields the desired result.   \endproof

Let $r(z)$ be as in Proposition \ref{refinement}. Let $c\in N^{|z|}$ be such that $N\models r(c)$. Then $\theta(x,c)$ does not divide over $M$ (since that is a property of $r(z)$ and does not depend on which realisation we are using). It follows by Proposition \ref{secondapprox} that there exist $\tau(z)\in r(z)$ and a finite set $A_{\tau}\subseteq N^{|x|}$ such that, for all $c^\prime\in \tau(M)$, there is some $a\in A_{\tau}$ such that $N\models \theta(a,c^\prime)$.

Let $j=|A_\tau|$. By Proposition \ref{refinement} and compactness, there is some $\psi(y)\in tp(b/M)$ such that, for all $B\subseteq \psi(M)$ with $|B|\leq j$, there exists $c^\prime\in\tau(M)$ such that $N\models \forall x(\theta(x,c^\prime)\rightarrow\varphi(x,d))$ for all $d\in B$.  

We complete the proof of Theorem \ref{mainresult} by showing that there is some $a\in A_\tau$ such that $N\models \varphi(a,b^\prime)$ for all $b^\prime\in \psi(M)$. Suppose not. Enumerate the elements of $A_\tau$ as $a_1,...,a_j$. Then, for each $i\in \{1,...,j\}$, there is some $b_i\in \psi(M)$ such that $N\models\neg\varphi(a_i,b_i)$. Then $B=\{b_1,...,b_j\}$ is a subset of $\psi(M)$ such that $|B|\leq j$. So there is some $c^\prime\in\tau(M)$ such that $N\models \forall x(\theta(x,c^\prime)\rightarrow\varphi(x,b_i))$ for all $i\in\{1,...,j\}$. Also there is some $i\in \{1,...,j\}$ such that $N\models \theta(a_i,c^\prime)$. Therefore there is some $i\in \{1,...,j\}$ such that $N\models \varphi(a_i,b_k)$ for all $k\in\{1,...,j\}$. This is a contradiction and the proof is finished.

\bibliographystyle{plain}

\bibliography{referencesBKup}

\begin{thebibliography}{1}

\bibitem{CK}
Artem Chernikov and Itay Kaplan.
\newblock Forking and dividing in {${\rm NTP}\sb 2$} theories.
\newblock {\em J. Symbolic Logic}, 77(1):1--20, 2012.

\bibitem{CS}
Artem Chernikov and Pierre Simon.
\newblock Externally definable sets and dependent pairs ii.
\newblock {\em Transactions of the American Mathematical Society}, 367:5217 --
  5235.

\bibitem{Shelah}
Saharon Shelah.
\newblock Dependent first order theories, continued.
\newblock {\em Israel Journal of Mathematics}, 173(1):1--60, 2009.

\bibitem{Sinv}
Pierre Simon.
\newblock Invariant types in nip theories.
\newblock {\em Journal of Mathematical Logic}, 0(0):1550006, 0.

\bibitem{Sdist}
Pierre Simon.
\newblock Distal and non-distal {NIP} theories.
\newblock {\em Ann. Pure Appl. Logic}, 164(3):294--318, 2013.

\bibitem{Sdp}
Pierre Simon.
\newblock Dp-minimality: Invariant types and dp-rank.
\newblock {\em The Journal of Symbolic Logic}, 79:1025--1045, 12 2014.

\bibitem{Sbook}
Pierre Simon.
\newblock {\em A Guide to NIP Theories}.
\newblock Cambridge University Press, 2015.
\newblock Cambridge Books Online.

\bibitem{SS}
Pierre Simon and Sergei Starchenko.
\newblock On forking and definability of types in some dp-minimal theories.
\newblock {\em The Journal of Symbolic Logic}, 79:1020--1024, 12 2014.

\end{thebibliography}

\end{document}